\documentclass[reqno]{amsart}
\usepackage{amsthm}
\usepackage{times,amssymb,amsmath,amsfonts,bbm,mathrsfs}
\usepackage[mathscr]{eucal}
\newcommand\nc\newcommand
\nc\bfa{{\boldsymbol a}}\nc\bfA{{\boldsymbol A}}\nc\cA{{\mathscr A}}
\nc\bfb{{\boldsymbol b}}\nc\bfB{{\boldsymbol B}}\nc\cB{{\mathscr B}}
\nc\bfc{{\boldsymbol c}}\nc\bfC{{\boldsymbol C}}\nc\cC{{\mathscr C}}
\nc\bfd{{\boldsymbol d}}\nc\bfD{{\boldsymbol D}}\nc\cD{{\mathscr D}}
\nc\bfe{{\boldsymbol e}}\nc\bfE{{\boldsymbol E}}\nc\cE{{\mathscr E}}
\nc\bff{{\boldsymbol f}}\nc\bfF{{\boldsymbol F}}\nc\cF{{\mathscr F}}
\nc\bfg{{\boldsymbol g}}\nc\bfG{{\boldsymbol G}}\nc\cG{{\mathscr G}}
\nc\bfh{{\boldsymbol h}}\nc\bfH{{\boldsymbol H}}\nc\cH{{\mathscr H}}
\nc\bfi{{\boldsymbol i}}\nc\bfI{{\boldsymbol I}}\nc\cI{{\mathcal I}}
\nc\bfj{{\boldsymbol j}}\nc\bfJ{{\boldsymbol J}}\nc\cJ{{\mathscr J}}
\nc\bfk{{\boldsymbol k}}\nc\bfK{{\boldsymbol K}}\nc\cK{{\mathscr K}}
\nc\bfl{{\boldsymbol l}}\nc\bfL{{\boldsymbol L}}\nc\cL{{\mathscr L}}
\nc\bfm{{\boldsymbol m}}\nc\bfM{{\boldsymbol M}}\nc{\cM}{{\mathscr M}}
\nc\bfn{{\boldsymbol n}}\nc\bfN{{\boldsymbol N}}\nc\cN{{\mathscr N}}
\nc\bfo{{\boldsymbol o}}\nc\bfO{{\boldsymbol O}}\nc\cO{{\mathscr O}}
\nc\bfp{{\boldsymbol p}}\nc\bfP{{\boldsymbol P}}\nc\cP{{\mathscr P}}
\nc\bfq{{\boldsymbol q}}\nc\bfQ{{\boldsymbol Q}}\nc\cQ{{\mathscr Q}}
\nc\bfr{{\boldsymbol r}}\nc\bfR{{\boldsymbol R}}\nc\cR{{\mathscr R}}
\nc\bfs{{\boldsymbol s}}\nc\bfS{{\boldsymbol S}}\nc\cS{{\mathscr S}}
\nc\bft{{\boldsymbol t}}\nc\bfT{{\boldsymbol T}}\nc\cT{{\mathscr T}}
\nc\bfu{{\boldsymbol u}}\nc\bfU{{\boldsymbol U}}\nc\cU{{\mathscr U}}
\nc\bfv{{\boldsymbol v}}\nc\bfV{{\boldsymbol V}}\nc\cV{{\mathscr V}}
\nc\bfw{{\boldsymbol w}}\nc\bfW{{\boldsymbol W}}\nc\cW{{\mathscr W}}
\nc\bfx{{\boldsymbol x}}\nc\bfX{{\boldsymbol X}}\nc\cX{{\mathscr X}}
\nc\bfy{{\boldsymbol y}}\nc\bfY{{\boldsymbol Y}}\nc\cY{{\mathscr Y}}
\nc\bfz{{\boldsymbol z}}\nc\bfZ{{\boldsymbol Z}}\nc\cZ{{\mathscr Z}}
\nc{\bb}{{\mathbbm{1}}}
\nc\RR{{\mathbb R}}
\nc\EE{{\mathbb E}}
\def\dd{\mathrm{d}}

\theoremstyle{remark}

\newcommand{\TT}{{\mathbb T}}
\newcommand{\ZZ}{{\mathbb Z}}
\newcommand{\KK}{{\mathbb K}}

\begin{document}

\title[Discrepancies and their means] 
{Discrepancies and their means}
\author[M.M. SKRIGANOV]{M.M. SKRIGANOV}

\address{St. Petersburg Department of the Steklov Mathematical Institute 
	of the Russian Academy of Sciences, 
	27, Fontanka, St.Petersburg 191023, Russia}

\email{maksim88138813@mail.ru}



 \begin{abstract}
It is shown that the discrepancy function for point distributions on 
a torus is expressed by an explicit formula in terms of its mean values on 
sub-tori.  
As an application of this formula, a simple proof of a theorem of 
Lev \cite{Lev} on the equivalence of  
$L_{\infty}$-- and  shifted $L_q$--discrepancies is given. 	
 \end{abstract}
 \maketitle

The point distribution problem on the 
$d$-dimensional 
torus $\TT^d =\RR^d/\ZZ^d$
is conveniently considered as a periodic problem on the covering space $\RR^d$.
For $X=(x_1,\dots ,x_d)\in\RR^d $ and $Y=(y_1,\dots ,y_d)\in\KK^d =[0,1]^d$,
we define the periodic discrepancy function by
\begin{align}
L(X,Y)=\chi(X,Y)-v(Y),
\label{eq:1}
\end{align}
where
$
\chi(X,Y)=\prod\nolimits_{j=1}^d\chi(x_j,y_j),\,\,
v(Y)=\prod\nolimits_{j=1}^d v(y_j), v(y)=y,
$
and
\begin{align}
\chi(x,y)=\begin{cases} 
	1,&\text{if $\{x\}<y$,}\\
	0,&\text{otherwise,}
\end{cases}
\label{eq:2}
\end{align}
here $\{x\}$ is the fractional part of $x\in\RR$, and $y\in [0,1]$.
It is clear that $\chi(X,Y)$ is the indicator function of the periodic 
collection of rectangular boxes 
\begin{align*}
\cB (Y)=\bigcup\nolimits_{(m_1,\dots, m_d)\in\ZZ^d}\prod\nolimits^d_{j=1} 
[m_j,y_j+m_j).
	\end{align*}

The mean value of the discrepancy function (1) has the form
\begin{align}
\lambda (X)=\int\nolimits_{\KK^d} L(X,Y)\,\dd Y =
\prod\nolimits_{j=1}^d (1-\{x_j\}) -2^{-d}.
\label{eq:4}
\end{align}

Let $[d]=(1,\dots,d)$ denote the set of coordinate indexes. For  subsets 
$J\subseteq [d]$ we introduce the partial discrepancy functions by
\begin{align}
L_J(X,Y)=\chi_J(X,Y)-v_J(Y),
\label{eq:5}
\end{align}
where
$
\chi_J(X,Y)=\prod\nolimits_{j\in J}\chi(x_j,y_j),\,\,
v_J(Y)=\prod\nolimits_{j\in J} v(y_j),
$
The corresponding mean values have the form 
\begin{align}
\lambda_J (X)=\int\nolimits_{\TT^{|J|}} L_J(X,Y_J)\,\dd Y_J =
\prod\nolimits_{j\in J} (1-\{x_j\}) -2^{-|J|},
\label{eq:6}
\end{align}
where $\TT^{|J|}\subseteq \TT^d$ denotes the sub-torus corresponding to the 
subset $J\subseteq [d]$, and $|J|$ denotes the number of elements of $J$.

The quantities \eqref{eq:5}, \eqref{eq:6} depend on the projections
$X_J=(x_j)_{j\in J}\in\RR^{|J|}, Y_J=(y_j)_{j\in J}\in\KK^{|J|}$, but not on 
the additional variables 
$X_{J'}=(x_j)_{j\in J'}\in\RR^{|J'|}, Y_{J'}=(y_j)_{j\in J'}\in\KK^{|J'|}$, 
where $J'=[d]\setminus J$ denotes the compliment of $J$. For $J=[d]$, we write
$L_{[d]}(X,Y)=L(X,Y),\, \lambda_{[d]} (X)=\lambda (X)$, and for the empty set, 
we put 
$L_{\emptyset}(X,Y)=0,\,  \lambda_{\emptyset} (X)=0$. 

{\bfseries Definition.} 
For a periodic function $f(X)=f(X_J,X_{J'}), X\in\RR^d$, and 
a vector 
$Y_J=(y_j)_{j\in J}\in \KK^{|J|}$, the {\it alternant} is defined by 
\begin{align}
f^{(alt)} (X\,|\,Y_J)=f^{(alt)} 
(X_J,X_{J'}\,|\,Y_J)=\sum\nolimits_{\Theta_J}(-1)^{|\Theta_J|} 
f(X_J-\Theta_J\cdot Y_J,X_{J'}),
\label{eq:7}
\end{align}
where $\Theta_J=(\theta_j)_{j\in J}\in \{0,1\}^{|J|}$ are the 
vertexes of the cube $\KK^{|J|}$, and summation in \eqref{eq:7} is taken over 
all such 
vertexes, $\Theta_J\cdot Y_J=(\theta_j y_j)_{j\in J}$ and $|\Theta_J |
=\sum\nolimits_{j\in J}\theta_j$.


{\bfseries Remark.} If $f(X)=\prod\nolimits_{j\in J}f_j(x_j)$, 
then
\begin{align}
f^{(alt)}(X\,|\,Y_J)=\prod\nolimits_{j\in J}f^{(alt)}_j(x_j\,|\,y_j)
=\prod\nolimits_{j\in J}(f(x_j)-f(x_j -y_j)) ,
\label{eq:8}
\end{align}
and if at least one of the functions $f_j$ is constant, then 
$f^{(alt)}=0$.

{\bfseries Theorem (Main Identity).} {\it The discrepancy function $L (X,Y)$ 
satisfies the 
identity}
\begin{align}
L (X,Y)=\sum\nolimits_{J\subseteq [d]}\, v_{J'}(Y) 
\,\lambda^{(alt)}_J(X\,|\,Y_J).
\label{eq:9}
\end{align}

\begin{proof}
	Let	$\omega(x)=\frac12 -\{x\}$, then $ \int\nolimits^1_0\omega(x)\dd x=0$.  
	We put
	\begin{align}
	\omega_J (X)=\prod\nolimits_{j\in J}\omega 
	(x_j)\qquad\text{and} \qquad\omega_{\emptyset} 
	(X)=0.
	\label{eq:10} 
	\end{align}
	
	The indicator function \eqref{eq:2} can be written in the form
	\begin{align}
	\begin{split}
		\chi (x,y)&=y-\{x\}+\{x-y\}\\&=y+\omega (x) -\omega 
		(x-y)=y+\omega^{(alt)}(x\,|\,y),	
	\end{split}
\label{eq:11}	
	\end{align}
	This formula can be proved by considering the graph of the function 
	$\{x\}-\{x-y\}, x\in\RR$.
	
	Substituting \eqref{eq:11} into \eqref{eq:1} and using \eqref{eq:10} and 
	\eqref{eq:8}, we obtain 
	\begin{align*}
	L (X,Y)=\sum\nolimits_{J\subseteq [d]}\, v_{J'}(Y) 
	\,\omega^{(alt)}_J(X\,|\,Y_J).
\end{align*}
	For the mean value \eqref{eq:6}, we find
	\begin{align}
	\lambda_J (X)=
	\prod\nolimits_{j\in J} (2^{-1} +\omega (x_j)) -2^{-|J|}=
	\sum\nolimits_{I\subseteq J}\, 2^{-|J\setminus I|} \,\,\omega_I (X).
	\label{eq:13}
	\end{align}

Let us calculate the alternant of the mean value \eqref{eq:13}. We have
\begin{align*}
\lambda^{(alt)}_J (X\,|\,Y_J)=
\omega^{(alt)}_J (X)+\sum\nolimits_{I\subset J}\, 2^{-|J\setminus I|} 
\,\,\omega^{(alt)}_I (X\,|\,Y_J).
\end{align*}
By the above Remark $\omega^{(alt)}_I (X\,|\,Y_J)=0$ for proper subsets 
$I\subset J$, since  
$\omega_J (X)=\prod\nolimits_{j\in I}\omega 
(x_j)\prod\nolimits_{j\in J\setminus I} 1$. 
Therefore, $\lambda^{(alt)}_J (X\,|\,Y_J)=\omega^{(alt)}_J(X\,|\,Y_J)$, and  
Theorem follows.
\end{proof}

We consider periodic point distributions $\cD$ on $\RR^d, \cD +M=\cD, 
M\in\ZZ^d,$ 
with a finite set of residues $\cD /\ZZ$. Notice that instead of point 
distributions,  
arbitrary periodic complex Borel measures on  $\RR^d$ finite on $\TT^d$ could 
be considered but
we do not consider such a generalization in order not to complicate the 
notation.

We define the local discrepancy
\begin{align}
L [\cD,Y]=\sum\nolimits_{X\in\cD /\ZZ^d} L(X,Y),
\label{eq:15}
\end{align}
and the $L_q$-- discrepancies 
\begin{align*}
L_q[\cD]=\left(\int\nolimits_{\KK^d} |L[\cD,Y]|^q\,\dd Y\right)^{1/q} ,\, 
0<q<\infty ,\quad\quad
L_{\infty}[\cD]=\sup\nolimits_{Y\in\KK^d} |L[\cD,Y]| .
\end{align*}

We also introduce the shifted discrepancies
\begin{align*}
L_q^{*}[\cD]=\sup\nolimits_{Z\in\TT^d}L_q[\cD +Z], \quad
\quad
L_{\infty}^{*}[\cD]=\sup\nolimits_{Z\in\TT^d}L_{\infty}[\cD +Z],
\end{align*}
and their mean values
\begin{align*}
\lambda_J [\cD]=\sum\nolimits_{X\in\cD /\ZZ^d}\lambda_J (X),
\quad\quad \lambda^*_J [\cD]=\sup\nolimits_{Z\in\KK^d}\lambda_J [\cD +Z].
\end{align*}

The above Theorem immediately implies the following.

{\bfseries Lemma 1.}  {\it The local discrepancy $L [\cD,Y]$ satisfies the 
identity}
\begin{align}
L [\cD,Y]=\sum\nolimits_{J\subseteq [d]}\, v_{J'}(Y) 
\,\lambda^{(alt)}_J \,[\cD\,|\,Y_J],
\label{eq:20}	
\end{align}
{\it where $\lambda^{(alt)}_J[\cD\,|\,Y_J]=
\sum\nolimits_{X\in\cD /\ZZ^d}\, \lambda^{(alt)}_J(X\,|\,Y_J)$.}

{\it The discrepancy $L_{\infty} [\cD]$ satisfies the inequality}
\begin{align}
L_{\infty} [\cD]\le\sum\nolimits_{J\subseteq [d]} 
\,2^{|J|}\,\lambda^*_J[\cD].
\label{eq:21}		
\end{align}
\begin{proof}
The identity \eqref{eq:20} follows from \eqref{eq:9}. The definition 
\eqref{eq:7} implies
$
\lambda^{(alt)}_J[\cD\,|\,Y_J]\le 2^{|J|}\,\lambda^*_J[\cD], 	
	$
and the inequality \eqref{eq:21} follows, since $0\le v_{J'}(Y)\le 1$.
\end{proof}

The $L_q^{*}[\cD]$-discrepancies can be easily estimated by means from below.

{\bfseries Lemma 2.}  {\it For $1\le q\le\infty$ and any subset $J\subseteq 
[d]$, 
the discrepancy $L_q^{*}[\cD]$ 
satisfies the 
inequality}
\begin{align}
L_q^{*}[\cD]\ge 2^{d-|J|}\,\lambda^*_J [\cD].
\label{eq:22}
\end{align}	
\begin{proof} We have
\begin{align*}
L_q^{*}[\cD]\ge L_1^{*}[\cD]&=
\sup\nolimits_{Z\in\TT^d}\int\nolimits_{\KK^d}|L[\cD +Z, Y]|\,\dd Y \ge
\sup\nolimits_{Z\in\TT^d}\left|\int\nolimits_{\KK^d}L[\cD +Z, Y]\,\dd Y 
\right|\\ 
&\ge \sup\nolimits_{Z\in\TT^d}\left|\lambda[\cD +Z]\right|=\lambda^*[\cD].
	\end{align*}
For simplicity, we put $Z=(Z_J,Z_{J'})\in\TT^d, Z_J\in\TT^{|J|}, 
Z_{J'}\in\TT^{|J'|}$, and continue 
\begin{align*}
\lambda^*[\cD]&=	
\sup\nolimits_{Z_J}\sup\nolimits_{Z_{J'}}
|\lambda[\cD +Z]|\ge \sup\nolimits_{Z_J}
\int\nolimits_{\KK^{|J'|}}|\lambda[\cD +Z]|\,\dd Z_{J'}\\ 
&\ge\sup\nolimits_{Z_J}\left |\int\nolimits_{\KK^{|J'|}}\lambda[\cD 
+Z]\,\dd Z_{J'} \right|=
2^{-|J'|}\sup\nolimits_{Z_J} |\lambda[\cD +Z_J]|=
2^{d-|J|}\,\lambda^*_J [\cD],
	\end{align*}
that completes the proof.
\end{proof}

The next simple fact is well--known, see, for example, \cite[p. 4]{Lev}. For 
completeness, we give 
a short proof.

{\bfseries Lemma 3.}  {\it The $L_{\infty}$-- and $L_{\infty}^{*}$-- 
discrepancies 
are equivalent:} 
	\begin{align}
		L_{\infty}[\cD]\,\le \,L_{\infty}^{*}[\cD]\,\le\, 3^d L_{\infty}[\cD].
		\label{eq:19a}
\end{align} 
\begin{proof} The left inequality \eqref{eq:19a} is obvious. Let us prove 
the right. For $y\in [0,1],\, z\in [0,1)$, we introduce the notation 
\begin{align*}
	\delta_{y,z}=\begin{cases} 
		1,&\text{if \,$y+z\ge 1$,}\\
		0,&\text{otherwise,}
	\end{cases}
\end{align*}
The shifted indicator function, see 
\eqref{eq:2} and \eqref{eq:11}, can be written in the form 
\begin{align*}
	\chi (x+z,y)=\chi (x,y+{z})[1-\delta_{y,z}]+\chi 
	(x,y+{z}-1)\,\delta_{y,z}-\chi (x,{z})+\chi (x,1)\,\delta_{y,z},
\end{align*}
and similarly $v(y)=v(y+z)[1-\delta_{y,z}]+v(y+z-1)\,\delta_{y,z}-v (z)+v 
(1)\,\delta_{y,z}.$	Moreover, each of these formulas contains at most three 
non-zero terms.
Substituting these formulas into the definitions \eqref{eq:1} and  
\eqref{eq:15}, we find that the shifted discrepancy can be written as the sum  
$
	L [\cD +Z,Y]=\sum\nolimits_k c_k\,L[\cD, V_k],
	\label{eq:15a}
$
with some vectors $V_k=V_k(Y,Z)\in\KK^d$ and coefficients $c_k=c_k(Y,Z)$ equal 
either 
$\pm 1$ or $0$.
	Moreover, the sum contain at most $3^d$ non-zero terms. This implies the 
	right inequality \eqref{eq:19a}.
	\end{proof}

Lev \cite{Lev} established the equivalence of $L_{\infty}$-- and 
$L_q^{*}$-- discrepancies: 
$
L_{\infty}[\cD]\,\cong \,L_q^{*}[\cD]
$,
for all $q\ge 1$ with the implicit constants depending only on the 
dimension $d$. Another proof of the equivalence was given later by
Kolountzakis, we refer to the survey article \cite{Chen} by Chen for a  
detailed discussion of these issues. 

The equivalence of the $L_{\infty}$-- and 
$L_q^{*}$-- discrepancies can be easily 
derived from the foregoing statements.
We will formulate and prove 
the corresponding result in the following somewhat more general form.

{\bfseries Corollary (Lev's Equivalence).}  {\it For $0<q<\infty$, the 
$L_{\infty}$-- and 
$L_q^{*}$-- discrepancies are equivalent}:
\begin{align} 
3^{-d}\,L^*_q[\cD]\,\le \,  \,L_{\infty}[\cD]\,\le\, 
C_{d,q}\,L^*_q[\cD],
\label{eq:23}
\end{align}
{\it where the constant}
\begin{align}
C_{d,q}=\begin{cases} 
	\,(5/2)^{d},&\text{if $1\le q<\infty$,}\\
	\,(5/2)^{d/q}\,3^{d/q\, -d},&\text{if $0<q<1.$}
\end{cases}
\label{eq:24}
\end{align}

We did not seek to obtain the best constant $C_{d,q}$ in \eqref{eq:23}, 
however, we note that 
the formula \eqref{eq:24} correctly reflects the order of the constant for 
large and small $q$.
\begin{proof}[Proof of Corollary.] 
	The proof consists of three steps.
	
	
	{\it (i)} For $0<q <{\infty}$, the lower bound \eqref{eq:23} follows 
	from Lemma 3:
	$
	L^*_q[\cD]\,\le L_{\infty}^{*}[\cD]\,\le\, 3^d L_{\infty}[\cD],
	$
	since $L^*_q[\cD]$ is a non-decreasing function of $q>0$.
	
	{\it (ii)} For $1\le q<\infty$, 
	the upper bound \eqref{eq:23} follows from Lemma~1 and Lemma~2.
	Substituting \eqref{eq:22} into \eqref{eq:21}, we obtain 
\begin{align*}	
L_{\infty} [\cD]\le 2^{-d}\,\left(\sum\nolimits_{J\subseteq [d]} \,
\,2^{2|J|}\right)\,\, L_q^{*}[\cD]=(5/2)^d\,\, L_q^{*}[\cD].
\end{align*}

{\it (iii)} Finally, for $0< q<1$, the upper bound \eqref{eq:23}  follows from 
the logarithmic convexity of $L_q^{*}[\cD]$	as a function of $q>0$.
By the standard interpolation at points $q<1<p$, we obtain
\begin{align*}	
	L_1^{*}[\cD]\le (L_q^{*}[\cD])^{q\frac{p-1}{p-q}}\, 
	(L_p^{*}[\cD])^{p\frac{1-q}{p-q}}.
\end{align*}
This inequality takes the form 
\begin{align}	
L_1^{*}[\cD]\le (L_q^{*}[\cD])^{q}\, 
(L_{\infty}^{*}[\cD])^{1-q},
\label{eq:27a}
\end{align}
as $p\to\infty$.
	
The bounds  
$(5/2)^{-d}\,L_{\infty} [\cD]\le L_1^{*}[\cD]$ and $L^*_{\infty}[\cD]\le 3^d\, 
L_{\infty}[\cD],$
 are already established.
Substituting these bounds into \eqref{eq:27a}, we obtain the upper bound 
\eqref{eq:23} with the constant \eqref{eq:24}.

The proof of Corollary is completed.
\end{proof}

\end{document}